\documentclass[10pt,a4paper]{article}

\usepackage{amsmath,amsthm,amsfonts,url}
\usepackage{mathabx}
\usepackage[all]{xy}
\usepackage{xcolor}
\usepackage{pgf,tikz,graphicx}
\usepackage{fancyhdr}
\usepackage{authblk}

\numberwithin{equation}{section}

\theoremstyle{plain}
\newtheorem{theorem}{Theorem}[section]
\newtheorem{lemma}[theorem]{Lemma}
\newtheorem{proposition}[theorem]{Proposition}
\newtheorem{corollary}[theorem]{Corollary}

\theoremstyle{definition}
\newtheorem{definition}[theorem]{Definition}

\newtheorem{remark}[theorem]{Remark}

\newcommand{\Set}{\mathsf{Set}}

\newcommand{\fpMod}{\textsf{Mod}_{\textsf{fp}}}
\newcommand{\V}{\mathcal{V}}

\newcommand {\subjclass}[1]{\textbf{Mathematics Subject Classification.} #1.}
\newcommand {\keyword}[1]{\textbf{Keywords.} #1.}

\title{From coextensive varieties to the Gaeta topos}
\author{William Zuluaga}
\date{}

\begin{document}
\maketitle
\begin{abstract}
\noindent
In this paper, we show that in every coextensive variety $\mathcal{V}$, the assignment that maps each algebra to its set of central elements is both functorial and representable. Furthermore, we prove that the full subcategory of finitely presented algebras in $\mathcal{V}$ is coextensive. Finally, we establish that if $\mathcal{V}$ is additionally $(\vec{0}, \vec{1})$-dense, the Gaeta topos classifies central-free $\mathcal{V}$-models.

\noindent
\noindent
\keyword{Gaeta topos, central elements, coextensive varieties}
\end{abstract}

\section{Introduction}\label{Introduction}
A category is called \emph{extensive} if it has finite coproducts and the canonical functors $1\rightarrow \mathsf{C}/0$ and $\mathsf{C}/A\times \mathsf{C}/B\rightarrow \mathsf{C}/(A+B)$ are equivalences \cite{LX,CMZ2016}. Besides $\mathsf{Set}$, well-known examples of extensive categories include toposes, the category $\mathsf{Top}$ of topological spaces and continuous functions, the category of affine schemes, and the category of Priestley spaces. It is worth mentioning that the last two examples are equivalent to the dual categories of algebras, namely, the categories of commutative rings with unity and the category of bounded distributive lattices, respectively. In general, varieties of algebras whose dual category is extensive are called \emph{coextensive}. In \cite{Z2021}, coextensive varieties were studied using Vaggione's \emph{theory of central elements} (\cite{BV2013}, \cite{SV2009}, \cite{V1999}). Loosely speaking, this theory, rooted in universal algebra, facilitates the study of direct decompositions of the algebras within a variety by examining concrete elements of these algebras, specifically the so-called \emph{central elements}. These elements can be seen as a generalization of idempotent elements in commutative rings with unity and complemented elements in bounded distributive lattices. Coextensive varieties are of interest because, as noted by \cite{L2008} and \cite{M2021}, they provide an appropriate setting for developing algebraic geometry. Moreover, several varieties related to non-classical logics—such as Heyting algebras and MV-algebras—are coextensive, as are Church varieties \cite{CS2015}, known for their connection with the variety of lambda abstraction algebras. This suggests that developing general techniques for coextensive varieties may be useful in several branches of mathematics.

Small extensive categories admit a particular subcanonical topology called the Gaeta topology, seemingly named in honor of the Spanish algebraic geometer Federico Gaeta. This topology is related to all possible decompositions of objects into finite coproducts, and the topos it generates is called the \emph{Gaeta topos}. In concrete examples (\cite{Za2017}, \cite{M2021}), it has been proved that the Gaeta topos classifies directly indecomposable algebras. In those instances, the theory whose models are precisely indecomposable algebras can be described in a straightforward manner. However, it is important to note that in the presentation of that theory, the set of equations that define central elements is always involved; thus, directly indecomposable algebras are those which cannot be decomposed by central elements. It is worth stressing that in arbitrary varieties, central elements do not have to be equationally defined (see Theorem 5 of \cite{V1996}). Therefore, considering the characterization given in \cite{Z2021} for coextensive varieties, it seems reasonable to ask: (1) Is it possible to propose a theory within a topos that precisely encompasses the notion of algebras that cannot be decomposed by central elements, namely, central-free models? And (2) Can we determine, using the tools provided by universal algebra and the theory of central elements, a suitable framework in which the Gaeta topos classifies central-free models? In this paper, we provide an affirmative answer to both of these questions.

This paper is organized as follows. Section \ref{Preliminaries} presents most of the definitions and basic results required for reading this work. Section \ref{Representability of the functor Z} is devoted to showing that in coextensive varieties, there is a natural notion of central-elements functor that coincides with the general results about coextensive categories obtained in \cite{G1998} and \cite{CPR2001}. Section \ref{The Gaeta topos and fp-coextensive varieties} is committed to show that the Gaeta topos classifies central-free models of $(\vec{0},\vec{1})$-dense coextensive varieties. 

The reader is assumed to be familiar with some standard topos theory and the basics of the internal logic of toposes as presented in \cite{MM2012} and \cite{J2002}. For usual notions in universal algebra, the reader may consult \cite{MMT1987}.

\section{Preliminaries}\label{Preliminaries}

\subsection{Notation and basic results}\label{Notation and basic results}

Let $A$ be a set and $N$ be a natural number. We write $\vec{a}$ for an element $(a_{1},\dots,a_{N})\in A^{N}$. If $f:A\rightarrow B$ is a function and $\vec{a}\in A^{N}$, then we write $f(\vec{a})$ for the element $(f(a_{1}),\dots,f(a_{N}))\in B^{N}$. If $X\subseteq A$ we write $f|_{X}$ for the restriction of $f$ to $X$, $\mathcal{P}(X)$ for the power set of $X$ and $f[X]$ for the image of $X$ through $f$. If $\vec{a} \in A^{N}$ and $\vec{b} \in B^{N}$, we write $[\vec{a}, \vec{b}]$ for the $N$-tuple $((a_{1} , b_{1} ), \dots, (a_{N} , b_{N})) \in (A \times B)^{N}$. If $g:A\times B \rightarrow C$ is a function and $[\vec{a}, \vec{b}]\in (A \times B)^{N}$ then we write $g(\vec{a},\vec{b})$ for the element $(g(a_{1},b_{1}),\dots,g(a_{N},b_{N}))\in C^{N}$.  If $\mathbf{A}$ is an algebra of a given type we denote its universe by $A$. If $X\subseteq A$ we denote by $\mathsf{Sg}^{\mathbf{A}}(X)$ the subalgebra of $\mathbf{A}$ generated by $X$. We also write $\mathsf{Con}(\mathbf{A})$ to denote the congruence lattice of $\mathbf{A}$. If $\theta \in \mathsf{Con}(\mathbf{A})$, and $\vec{a}\in A^{N}$ we write $\vec{a}/\theta$ for the $N$-tuple $(a_{1}/\theta , \dots, a_{N}/\theta)\in (A/\theta)^{N}$. The universal congruence on $\mathbf{A}$ and identity congruence on $\mathbf{A}$ are denoted by $\nabla^{\mathbf{A}}$ and $\Delta^{\mathbf{A}}$, respectively. If $S\subseteq A\times A$, we write $\mathsf{Cg}^{\mathbf{A}}(S)$ for the congruence generated by $S$. We also write $\mathsf{Cg}^{\mathbf{A}}(\vec{a},\vec{b})$, for the congruence generated by all pairs $(a_{1} , b_{1} ), \dots, (a_{N} , b_{N})$ where $\vec{a}, \vec{b}\in A^{N}$.  A pair $\theta, \delta\in \mathsf{Con}(\mathbf{A})$ which satisfies $\theta\cap\delta = \Delta^{\mathbf{A}}$ and $\theta\circ\delta = \nabla^{\mathbf{A}}$ is called a pair of \emph{complementary factor congruences} (where $\circ$ denotes the usual composition of congruences). We write $\theta \diamond \delta$ in $\mathsf{Con}(\mathbf{A})$ to denote that $\theta$ and $\delta$ are complementary factor congruences of $\mathbf{A}$. We say that $\theta\in \mathsf{Con}(\mathbf{A})$ is a \emph{factor congruence} if there exists a $\delta\in \mathsf{Con}(\mathbf{A})$ such that $\theta \diamond \delta$. It is well known that factor congruences concentrate all the information about decompositions of algebras into direct products. We use $\mathsf{FC}(\mathbf{A})$ to denote the set of factor congruences of $\mathbf{A}$. A variety $\V$ has the \emph{Fraser-Horn property}  \cite{FH1970} if for every $\mathbf{A}_1,\mathbf{A}_2 \in \V$, every congruence $\theta$ in $\mathbf{A}_1\times\mathbf{A}_2 $ is the product congruence $\theta_1 \times \theta_2$ for some congruences $\theta_1$ of $\mathbf{A}_1$ and $\theta_2$ of $\mathbf{A}_2$. If $\theta,\lambda \in \mathsf{Con}(\mathbf{A})$ and $\theta \subseteq \lambda$, we write $\lambda/\theta$ for the set of pairs $(x/\theta,y/\theta)$ of $A/\theta$ such that $(x,y)\in \lambda$. For any homomorphism $g:\mathbf{A}\rightarrow \mathbf{B}$ we write  $\mathsf{Ker}(g)$ for the kernel of $g$. If $\mathbf{A}$ is an algebra of type $\mathcal{F}=\{f_{1},\ldots ,f_{m}\}$, when required, we will write its type as an $m$-tuple $(c_{1},\ldots, c_{m})$, where $c_{j}$ denotes the arity of $f_{j}$, with $1\leq j\leq m$.

The following result (probably folklore), provides a description of the factor congruences of the quotients of an arbitrary algebra. It is is a straightforward consequence of Theorems 7.5, 6.15 and 6.20 of \cite{BS1981}. 

\begin{lemma}\label{factor congruence quotinents}
Let $\mathbf{A}$ be an algebra and let $\theta\in \mathsf{Con}(\mathbf{A})$. Consider the sets \[P_{\theta}=\{(\lambda,\mu)\mid\theta\subseteq\lambda,\mu;\;\lambda\cap\mu=\theta;\;\lambda\circ\mu=\nabla^{\mathbf{A}}\}\] 
and 
\[Z_{\theta}=\{(\alpha,\beta) \in \mathsf{FC}(\mathbf{A}/\theta)^{2}\colon \alpha \diamond \beta\}.\]
Then, the assignment $(\lambda,\mu)\mapsto (\lambda/\theta,\mu/\theta)$ defines a bijection between $P_{\theta}$ and $Z_{\theta}$. 
\end{lemma}

Given a set $X$ of variables and a variety $\mathcal{V}$, we use $\mathbf{T}_{\V}(X)$ for the term algebra of $\V$ over $X$ and $\mathbf{F}_{\mathcal{V}}(X)$ for the free algebra of $\mathcal{V}$ freely generated by $X$. Essentially, the free algebra $\mathbf{F}_{\mathcal{V}}(X)$ is built as a quotient of the term algebra $\mathbf{T}_{\V}(X)$ and its elements are congruence classes of terms equivalent in $\V$. If there is no place to confusion, we identify the universe of $\mathbf{F}_{\V}(X)$ with a set
of its representatives, i.e., with a set of terms in variables $X$. If $X=\{x_{1},\dots,x_{m}\}$ with $m$ a non-negative integer and if no clarification is needed, then we write $\mathbf{T}_{\V}(m)$ and $\mathbf{F}_{\V}(m)$ instead of $\mathbf{T}_{\V}(\{x_{1},\dots,x_{m}\})$ and $\mathbf{F}_{\V}(\{x_{1},\dots,x_{m}\})$, respectively. We recall that an algebra $\mathbf{A}$ in $\V$ is a \emph{finitely generated free algebra} if it is isomorphic to $\mathbf{F}_{\V}(m)$ for some finite $m$, and \emph{finitely presented} if it is isomorphic to an algebra of the form $\mathbf{F}_{\V}(k)/\theta$, for some $k$ finite and $\theta$ finitely generated congruence on $\mathbf{F}_{\V}(k)$.

The following Lemma is a key result that we will employ repeatedly along Section \ref{The Gaeta topos and fp-coextensive varieties}. For its proof, the reader may consult \cite{V1996_1}.

\begin{lemma}\label{vey useful lemma}
Let $\V$ be a variety and let $X$ be a set of variables. Let $r, r_1 , \ldots , r_m ,$ $s, s_1 , \ldots , s_m \in T_{\V}(X)$. Then, the following are equivalent:
\begin{itemize}
\item[(1)] $(r,s)\in \mathsf{Cg}^{\mathbf{F}_{\mathcal{V}}(X)}(\vec{r},\vec{s})$;
\item[(2)] $\V\models \vec{r}=\vec{s}\rightarrow r=s$.
\end{itemize}

\end{lemma}

Let $\mathcal{L}$ be a first order language. If $\mathcal{K}$ is a class of $\mathcal{L}$-structures and $R\in \mathcal{L}$ is a $n$-ary relation symbol, we say that a formula $\varphi (x_{1},\dots,x_{n})$ \emph{defines} $R$ \emph{in} $\mathcal{K}$ if
\begin{equation*}
\mathcal{K}\vDash \varphi (\vec{x})\leftrightarrow R(\vec{x})\text{.}
\end{equation*}
In particular, if a $\mathcal{L}$-formula $\varphi (\vec{x})$ which is the conjunction of a finite number of equations in the variables $x_{1},\dots,x_{n}$ defines $R$, we say that $R$ is \emph{equationally definable}.

Finally, we stress that all varieties in this paper are assumed to have at least a constant symbol in their language.

\subsection{Central Elements}\label{Central Elements}

By \emph{a variety with $\vec{0}$ and $\vec{1}$} we understand a variety $\V$ in which there are $0$-ary terms $0_{1}$, $\ldots$ , $0_{N}$, $1_{1}$, $\ldots$ , $1_{N}$ such that $\V \models \vec{0}= \vec{1}\rightarrow x= y$, where $\vec{0}=(0_{1}, ..., 0_{N})$ and $\vec{1}=(1_{1}, ..., 1_{N})$.  Varieties with $\vec{0}$ and $\vec{1}$ are exactly the ones in which the trivial algebra embeds only in itself. They can also be described, as noted by J. Koll\'ar \cite{K1979}, as the
varieties in which the universal congruences $\nabla^{\mathbf{A}}$ are finitely generated for every algebra $\mathbf{A}$ in the variety. 

If $\mathbf{A}\in \V$ then we say that $\vec{e}=(e_{1}, ..., e_{N})\in A^{N}$ is a \emph{central element} of $\mathbf{A}$ if there exists an isomorphism	$\tau: \mathbf{A}\rightarrow \mathbf{A}_{1}\times \mathbf{A}_{2}$, such that $\tau(\vec{e})=[\vec{0}_{\mathbf{A}_{1}}, \vec{1}_{\mathbf{A}_{2}}]$. Also, we say that $\vec{e}$ and $\vec{f}$ are a \emph{pair of complementary central elements} of $\mathbf{A}$ if there exists an isomorphism $\tau: \mathbf{A}\rightarrow \mathbf{A}_{1}\times \mathbf{A}_{2}$ such that $\tau(\vec{e})=[\vec{0}_{\mathbf{A}_{1}}, \vec{1}_{\mathbf{A}_{2}}]$ and $\tau(\vec{f})=[\vec{1}_{\mathbf{A}_{1}}, \vec{0}_{\mathbf{A}_{2}}]$. We write $Z(\mathbf{A})$ to denote the set of central elements of $A$ and $\vec{e}\diamond_{\mathbf{A}} \vec{f}$ to denote that $\vec{e}$ and $\vec{f}$ are complementary central elements of $\mathbf{A}$. We say that a variety $\V$ with $\vec{0}$ and $\vec{1}$ has \emph{Boolean Factor Congruences} (BFC) if the set of factor congruences of any algebra of $\V$ is a Boolean sublattice of its congruence lattice. We stress that such assumptions are enough for making that central elements concentrate all the information about direct-product decompositions of the members of a given variety. This is a generalization of what happens with idempotent elements and complemented elements in commutative rings with unit and bounded distributive lattices, respectively. Such an assertion on central elements is justified by Theorem 1.1 of \cite{SV2009}. It states that in every variety $\V$ with $\vec{0}$ and $\vec{1}$ and BFC the assignment 
\begin{displaymath}
\begin{array}{cccc}
\varphi: & Z(\mathbf{A}) & \longrightarrow & \mathsf{FC}(\mathbf{A})
\\
 & e & \longmapsto & \theta,
\end{array}
\end{displaymath}

where $\theta$ is the unique factor congruence of $\mathbf{A}$ such that $[\vec{e}, \vec{0}] \in \theta$ and $[\vec{e}, \vec{1}]\in \theta^{\ast}$, with $\theta\diamond \theta^{\ast}$, is a bijection, for every $\mathbf{A}\in \V$. If $\vec{e}\in Z(\mathbf{A})$, we write $\theta_{\vec{0}, \vec{e}}^{\mathbf{A}}$ for $\varphi(\vec{e})$ and $\theta_{\vec{1}, \vec{e}}^{\mathbf{A}}$ for $\varphi(\vec{e})^{\ast}$. Moreover, such a bijection, allows to endow $Z(\mathbf{A})$ with some operations in such a way that the algebra  $\textbf{Z}(\mathbf{A})=(Z(\mathbf{A}),\wedge_{\mathbf{A}},\vee_{\mathbf{A}}, ^{c_{\mathbf{A}}},\vec{0}^{\mathbf{A}},\vec{1}^{\mathbf{A}})$ becomes a Boolean algebra which is isomorphic to the Boolean algebra $(\mathsf{FC}(\mathbf{A}), \vee, \cap, ^{\ast},\Delta^{\mathbf{A}},\nabla^{\mathbf{A}})$. These operations are defined as follows: given $\vec{e}\in Z(\mathbf{A})$, the \emph{complement $\vec{e}^{c_{\mathbf{A}}}$} of $\vec{e}$, is the only solution to the equations $[\vec{z}, \vec{1}]\in \theta^{\mathbf{A}}_{\vec{0},\vec{e}}$ and $[\vec{z}, \vec{0}]\in \theta^{\mathbf{A}}_{\vec{1},\vec{e}}$. Given $\vec{e}, \vec{f}\in Z(\mathbf{A})$, the \emph{infimum} $\vec{e}\wedge_{\mathbf{A}}\vec{f}$ is the only solution to the equations $[\vec{z}, \vec{0}]\in \theta^{\mathbf{A}}_{\vec{0},\vec{e}}\cap \theta^{\mathbf{A}}_{\vec{0},\vec{f}}$ and $[\vec{z}, \vec{1}]\in \theta^{\mathbf{A}}_{\vec{1},\vec{e}}\vee \theta^{\mathbf{A}}_{\vec{1},\vec{f}}$. Finally, the \emph{supremum} $\vec{e}\vee_{\mathbf{A}}\vec{f}$ is the only solution to the equations $[\vec{z}, \vec{0}]\in \theta^{\mathbf{A}}_{\vec{0},\vec{e}}\vee \theta^{\mathbf{A}}_{\vec{0},\vec{f}}$ and $[\vec{z}, \vec{1}]\in \theta^{\mathbf{A}}_{\vec{1},\vec{e}}\cap \theta^{\mathbf{A}}_{\vec{1},\vec{f}}$.

Let $\V$ be a variety with BFC. If $\mathbf{A},\mathbf{B}\in \V$ and $f:\mathbf{A}\rightarrow \mathbf{B}$ is a homomorphism, we say that $f$ \emph{preserves central elements} if the map $f:Z(\mathbf{A})\rightarrow Z(\mathbf{B})$ is well defined; that is to say, for every $\vec{e}\in Z(\mathbf{A})$,  we have $f(\vec{e})\in Z(\mathbf{B})$. We say that $f$ \emph{preserves complementary central elements} if it preserves central elements and for every $\vec{e}_{1}, \vec{e}_{2}\in Z(\mathbf{A})$, 
\[\vec{e}_{1}\diamond_{\mathbf{A}}\vec{e}_{2} \Rightarrow f(\vec{e}_{1})\diamond_{\mathbf{B}}f(\vec{e}_{2}). \]

We say that a variety with BFC is \emph{center stable} if every homomorphism preserves central elements and we say that  it is \emph{stable by complements} if every homomorphism preserves complementary central elements. In \cite{Z2021} it was shown that these notions are not equivalent. 

\section{The central-elements functor}\label{Representability of the functor Z}

In this section we prove that in a coextensive variety $\V$, there is a concrete notion of a ``central-elements'' functor $Z:\V\rightarrow \Set$ which further, essentially coincides with the representable $\mathcal{V}(\mathbf{0}\times \mathbf{0},-)$. Our aim is to establish a connection between the results of \cite{G1998} and \cite{CPR2001} (which are given in an abstract categorical setting) and the theory of central elements, in such a way, that the latter may serve as an universal algebraic bridge between the general theory of coextensive varieties and the well-known cases studied along the literature. 

If $\V$ is a variety with BFC, the associated algebraic category will be also denoted by $\V$. The initial and terminal algebras of $\mathcal{V}$ will be denoted by $\mathbf{0}$ and $\mathbf{1}$, respectively. If $\mathbf{A}\in \mathcal{V}$ we write $\text{!`}_{\mathbf{A}}:\mathbf{0}\rightarrow \mathbf{A}$ for the unique morphism from $\mathbf{0}$ to $\mathbf{A}$ in $\mathcal{V}$. Further, if $\theta\in \mathsf{Con}(\mathbf{A})$ we write $\text{!`}_{\theta}$ for the unique morphism from $\mathbf{0}$ to $\mathbf{A}/\theta$. In particular, if $\vec{e}\in Z(\mathbf{A})$ we write $\text{!`}_{\vec{0},\vec{e}}$ and $\text{!`}_{\vec{1},\vec{e}}$ for the unique morphisms from  $\mathbf{0}$ to $\mathbf{A}/\mathsf{Cg}^{\mathbf{A}}(\vec{0},\vec{e})$ and from  $\mathbf{0}$ to $\mathbf{A}/\mathsf{Cg}^{\mathbf{A}}(\vec{1},\vec{e})$, respectively. Finally, we recall that due to the fact that $\V$ is assumed to have at least a constant symbol then $\mathbf{0}$ is isomorphic to $\mathbf{F}_{\V}(\emptyset)$.
\\

Let $\mathsf{C}$ be a category with binary sums (coproducts). We say that a morphism $Y\rightarrow X$ of $\mathsf{C}$ is a \emph{summand} of $X$ if there exists an arrow $X\leftarrow Z$ in $\mathsf{C}$ such that the diagram $Y\rightarrow X\leftarrow Z$ is a coproduct. We recall that a category $\mathsf{C}$ with finite sums is called \emph{extensive} if for each pair of objects $X$ and $Y$ of $\mathsf{C}$ the canonical functor $+: \mathsf{C}/X \times \mathsf{C}/Y \rightarrow \mathsf{C}/(X + Y)$ is an equivalence. In particular, if $\mathsf{C}$ has a terminal object $1$, from Proposition 4.1 of \cite{CW1993}, it follows that $\mathsf{C}$ is extensive if and only if the canonical functors $1\rightarrow \mathsf{C}/0$ and $\mathsf{C}/(1+1)\rightarrow \mathsf{C}\times \mathsf{C}$ are equivalences. In \cite{G1998} it was proved that in any extensive category $\mathsf{C}$ (without the need of a terminal object) the summands of an object $X$ in $\mathsf{C}$ form a Boolean algebra $X_{\ast}$. Moreover, if $f:X\rightarrow Y$ is an arrow in $\mathsf{C}$, then the pullback of a summand of $Y$ with $f$ induces a homomorphism of Boolean algebras $f_{\ast}:Y_{\ast}\rightarrow X_{\ast}$. This makes that the assignments $X\mapsto X_{\ast}$ and $f\mapsto f_{\ast}$ determine a product-preserving functor $(-)_{\ast}:\mathsf{C}^{\mathsf{op}}\rightarrow \mathsf{Boole}$. Let $\mathcal{O}$ be the forgetful functor from $\mathsf{Boole}$ to $\mathsf{Set}$. Observe that if $\mathsf{C}$ has a terminal object, it turns out that the composite $\mathsf{C}^{\mathsf{op}}\xrightarrow{(-)_{\ast}} \mathsf{Boole}\xrightarrow{\mathcal{O}} \mathsf{Set}$, is representable by $\mathsf{C}(-,1+1)$. We say that a category $\mathsf{D}$ is \emph{coextensive} if its opposite is extensive. In \cite{CPR2001}, it was shown that under the assumption of $\mathsf{D}$ is a coextensive locally presentable category, the functor $(-)_{\ast}$ possesses a left adjoint $(-)^{\ast}$ that preserves finite products. Furthermore, it was shown that the adjunction $(-)^{\ast}\dashv (-)_{\ast}$ is essentially unique. 
Having into account that every variety of algebras together with the homomorphisms between its members form a locally presentable category, it is clear that for coextensive varieties the previous results immediately take the following form: 

\begin{proposition}\label{prop: coextensive Gates+Carboni}
Let $\V$ be a coextensive variety and consider the representable functor $\mathcal{V}(\mathbf{0}\times \mathbf{0},-):\V\rightarrow \mathsf{Set}$. Then, the following hold:
\begin{enumerate}
\item There exists a functor $(-)_{\ast}:\mathcal{V}\rightarrow \mathsf{Boole}$ with a finite product-preserving left adjoint $(-)^{\ast}$. Moreover, such an adjunction is essentially unique.
\item The representable $\mathcal{V}(\mathbf{0}\times \mathbf{0},-)$ factors through $(-)_{\ast}$ and the forgetful functor from $\mathsf{Boole}$ to $\Set$. 
\end{enumerate}
\end{proposition}
We recall that a variety $\V$ is a \emph{Pierce variety} \cite{V1996} if there exist a positive natural number $N$, $0$-ary terms $0_1$, $\ldots$, $0_N$, $1_1$, $\ldots$, $1_N$ and a term $U (x, y, \vec{z}, \vec{w})$ such that the identities
\begin{displaymath}
\begin{array}{ccc}
U (x, y, \vec{0}, \vec{1})=x & \text{and} &  U (x, y, \vec{1}, \vec{0})=y,
\end{array}
\end{displaymath}
hold in $\V$. Such a term is called \emph{decomposition term} \cite{BV2013}. It is worth mentioning that in a Pierce variety $\V$, it is also true that $\theta^{\mathbf{A}}_{\vec{0},\vec{e}}=\mathsf{Cg}^{\mathbf{A}}(\vec{0},\vec{e})$, for every $\mathbf{A}\in \V$ and every $\vec{e}\in Z(\mathbf{A})$ (see \cite{BV2016} for details). Furthermore, from Proposition 3.2 of \cite{BV2013} we have that for every $\mathbf{A}\in \mathcal{V}$ and $\vec{e}\in \mathcal{V}$:
\begin{equation}\label{Ecucentrals}
\mathsf{Cg}^{\mathbf{A}}(\vec{e},\vec{0})=\{(a,b)\in A^2\colon U(a,b,\vec{e},\vec{1})=U(a,b,\vec{1},\vec{e})\}.
\end{equation}

Let $\V$ be a Pierce variety and let $\mathbf{A}\in \V$. If $\mathbf{A}$ is isomorphic to a product $\mathbf{A}_{0}\times \mathbf{A}_{1}$, there exists a unique pair $\theta,\delta \in \mathsf{FC}(\mathbf{A})$, with $\theta\diamond_{\mathbf{A}}\delta$, such that $\mathbf{A}_{0}$ is isomorphic to $\mathbf{A}/\theta$ and $\mathbf{A}_{1}$ is isomorphic to $\mathbf{A}/\delta$. The latter, by general reasons, implies that the homomorphism $\text{!`}_{\theta}\times \text{!`}_{\delta}:\mathbf{0}\times\mathbf{0}\rightarrow \mathbf{A}$ is the only arrow of $\V$ that makes both of the squares of the following diagram commute. 

\begin{displaymath}
\xymatrix{
\mathbf{0} \ar[d]_-{\text{!`}_{\theta}} & \ar[l]_-{\pi_{0}} \mathbf{0}\times\mathbf{0} \ar[d]^-{\text{!`}_{\theta}\times \text{!`}_{\delta}} \ar[r]^-{\pi_{1}}  & \mathbf{0} \ar[d]^-{\text{!`}_{\delta}}
\\
\mathbf{A}/\theta & \ar[r]_-{p_{1}} \ar[l]^-{p_{0}} \mathbf{A} & \mathbf{A}/\delta
}
\end{displaymath}

Considering that there is a bijection between $\mathsf{FC}(\mathbf{A})$ and the set of central elements of $\mathbf{A}$, then there exists a unique $\vec{e}\in \mathbf{A}$ such that $\theta=\mathsf{Cg}^{\mathbf{A}}(\vec{0},\vec{e})$ and $\delta=\mathsf{Cg}^{\mathbf{A}}(\vec{1},\vec{e})$. So, the arguments of above suggests that the assignment $\vec{e}\mapsto \text{!`}_{\theta}\times \text{!`}_{\delta}$ establishes a canonical map $\mu_{\mathbf{A}}:Z(\mathbf{A})\rightarrow \mathcal{V}(\mathbf{0}\times \mathbf{0}, \mathbf{A})$.

Let $f:\mathbf{A}\rightarrow \mathbf{B}$ be a homomorphism. Recall that in general, the restriction of $f$ to $Z(\mathbf{A})$ does not have $Z(\mathbf{B})$ as an image (see p.6 of \cite{Z2021}), so our next goal is to provide sufficient conditions on a variety $\V$ with BFC in order that: (1) the assignments $\mathbf{A}\mapsto Z(\mathbf{A})$ and $f\mapsto f|_{Z(\mathbf{A})}$ determine a functor $Z:\mathcal{V}\rightarrow \Set$, (2) the functions $\mu_{\mathbf{A}}$ determine a canonical natural transformation $\mu: Z\rightarrow \mathcal{V}(\mathbf{0}\times \mathbf{0}, -)$ and (3) there exists a natural transformation $\varphi: \mathcal{V}(\mathbf{0}\times \mathbf{0}, -)\rightarrow Z$.  Such a functor $Z$ will be called the \emph{central-elements} functor of $\mathcal{V}$.  


\begin{lemma}\label{Lemma suficiente}
Let $\mathcal{V}$ be a center stable variety with the Fraser-Horn property. Then the following hold:
\begin{enumerate}
\item For every $\mathbf{A},\mathbf{B}\in \V$ and homomorphism $f:\mathbf{A}\rightarrow \mathbf{B}$, the assignments $\mathbf{A}\mapsto Z(\mathbf{A})$ and $f \mapsto f|_{Z(\mathbf{A})}$ determine a functor $Z:\mathcal{V}\rightarrow \Set$.
\item There exist a pair of canonical maps $\mu:\mathcal{V}(\mathbf{0}\times \mathbf{0}, -)\rightarrow Z$ and $\varphi: Z\rightarrow \mathcal{V}(\mathbf{0}\times \mathbf{0}, -)$. 
\end{enumerate}
\end{lemma}
\begin{proof}
The first part follows from the definition of center stable variety. For the second part, notice that since $\V$ has the Fraser-Horn property, by Corollary 4 of \cite{V1999}, for every $\mathbf{A}\in \V$ and $\vec{e}\in Z(\mathbf{A})$, $\theta_{\vec{0},\vec{e}}^{\mathbf{A}}=\mathsf{Cg}^{\mathbf{A}}(\vec{0},\vec{e})$, so let $\mathbf{A}\in \mathcal{V}$ and consider the assignments $\mu_{\mathbf{A}}:Z(\mathbf{A}) \rightarrow \V(\mathbf{0}\times\mathbf{0}, \mathbf{A})$ and $\varphi_{\mathbf{A}}:\V(\mathbf{0}\times\mathbf{0}, \mathbf{A})\rightarrow Z(\mathbf{A})$, defined by $\mu_{\mathbf{A}}(\vec{e})=\text{!`}_{\vec{0},\vec{e}}\times \text{!`}_{\vec{1},\vec{e}}$ and $\varphi_{\mathbf{A}}(g)=g[\vec{0},\vec{1}]$, respectively. 
\begin{displaymath}
\xymatrix{
\mathbf{0} \ar[d]_-{\text{!`}_{\vec{0},\vec{e}}} & \ar[l]_-{\pi_{0}} \mathbf{0}\times\mathbf{0} \ar[d]^-{\text{!`}_{\vec{0},\vec{e}}\times \text{!`}_{\vec{1},\vec{e}}} \ar[r]^-{\pi_{1}}  & \mathbf{0} \ar[d]^-{\text{!`}_{\vec{1},\vec{e}}}
\\
\mathbf{A}/\mathsf{Cg}^{\mathbf{A}}(\vec{0},\vec{e}) & \ar[r]_-{p_{1}} \ar[l]^-{p_{0}} \mathbf{A} & \mathbf{A}/\mathsf{Cg}^{\mathbf{A}}(\vec{1},\vec{e})
}
\end{displaymath}
Since every $\vec{e}\in Z(\mathbf{A})$ induces a product decomposition of $\mathbf{A}$, then we get that $\text{!`}_{\vec{0},\vec{e}}\times \text{!`}_{\vec{1},\vec{e}}$ is a homomorphism in $\V(\mathbf{0}\times \mathbf{0}, \mathbf{A})$ making both squares of the diagram of above commute, thus $\mu_{\mathbf{A}}$ is also well defined. Similarly, because $[\vec{0},\vec{1}]\in Z(\mathbf{0}\times \mathbf{0})$, then the center stability of $\mathcal{V}$ makes that $g[\vec{0},\vec{1}]\in Z(\mathbf{A})$ for every $g\in \V(\mathbf{0}\times \mathbf{0}, \mathbf{A})$, so $\varphi_{\mathbf{A}}$ is well defined. The proof of the naturality of $\mu$ and $\varphi$ is straightforward.
\end{proof}
Now, we recall some facts about coextensive varieties which will required for proving the remaining results of this paper. It is well known that every variety (as a category) has all limits. Therefore, as a restricted dual of Propositions 2.2 and 4.1 of \cite{CW1993} we obtain the following result.

\begin{proposition}\label{coextensive with zero}
A variety $\mathcal{V}$ is coextensive if and only if it has pushouts along projections and every commutative diagram
\begin{displaymath}
\xymatrix{
\mathbf{0} \ar[d]_-{\text{!`}_{\mathbf{A}_{0}}} & \ar[l]_-{\pi_{0}} \ar[r]^-{\pi_{1}} \mathbf{0} \times \mathbf{0} \ar[d]_-{f} & \mathbf{0} \ar[d]^-{\text{!`}_{\mathbf{A}_{1}}} 
\\
\mathbf{A}_{0} & \ar[l]^-{g_{0}} \ar[r]_-{g_{1}} \mathbf{A}_{0}\times \mathbf{A}_{1} & \mathbf{A}_{1}
}
\end{displaymath}
comprises a pair of pushout squares in $\mathcal{V}$ just when the bottom row is a product diagram in $\mathcal{V}$.
\end{proposition}


\begin{theorem}[\cite{Z2021}]\label{charcoextensivity}
Let $\V$ be a variety. Then, the following are equivalent:
\begin{enumerate}
\item[(1)] $\mathcal{V}$ is coextensive.
\item[(2)] $\mathcal{V}$ is a Pierce variety in which the relation $\vec{e}\diamond_{\mathbf{A}} \vec{f}$ is equationally definable.  
\item[(3)] $\mathcal{V}$ is a Pierce variety stable by complements.
\end{enumerate}
\end{theorem}
Finally, we show that under the hypothesis that $\V$ is coextensive, the canonical $\mu$ of Lemma \ref{Lemma suficiente} is in fact a natural isomorphism.

\begin{theorem}\label{centrals bijection}
Let $\V$ be a coextensive variety. Then, for every $\mathbf{A}\in \mathcal{V}$ there is a bijection between $Z(\mathbf{A})$ and $\V(\mathbf{0}\times\mathbf{0}, \mathbf{A})$. Moreover, such a bijection is natural.
\end{theorem}
\begin{proof}
Observe that since $\V$ is coextensive, from Theorem \ref{charcoextensivity} (3) the hypotheses of Lemma \ref{Lemma suficiente} hold, so we can consider the canonical maps $\varphi: Z\rightarrow \mathcal{V}(\mathbf{0}\times \mathbf{0}, -)$ and $\mu:\mathcal{V}(\mathbf{0}\times \mathbf{0}, -)\rightarrow Z$ as defined in the aforementioned Lemma. We claim that under our assumptions, $\varphi$ and $\mu$ are inverse of each other.  In order to prove our claim, let $\mathbf{A}\in \V$, $\vec{e}\in Z(\mathbf{A})$ and let us write $h=\text{!`}_{\vec{0},\vec{e}}\times \text{!`}_{\vec{1},\vec{e}}$. It is clear that $\varphi_{\mathbf{A}}(\mu_{\mathbf{A}}(\vec{e}))=h[\vec{0},\vec{1}]$. 
\begin{displaymath}
\xymatrix{
\mathbf{0} \ar[d]_-{\text{!`}_{\vec{0},\vec{e}}} & \ar[l]_-{\pi_{0}} \mathbf{0}\times\mathbf{0} \ar[d]^-{h} \ar[r]^-{\pi_{1}}  & \mathbf{0} \ar[d]^-{\text{!`}_{\vec{1},\vec{e}}}
\\
\mathbf{P}_{0} & \ar[r]_-{j_{1}} \ar[l]^-{j_{0}} \mathbf{A} & \mathbf{P}_{1}
}
\end{displaymath}
Now, if $\mathbf{P}_{k}$ denote the pushouts of $\pi_{k}$ along $h$, with $1\leq k\leq 2$, then from Lemma 2.3 of \cite{Z2021}: 
\[
\begin{array}{ccc}
\mathbf{P}_{0}\cong \mathbf{A}/\mathsf{Cg}^{\mathbf{A}}(\vec{0},\varphi_{\mathbf{A}}(\mu_{\mathbf{A}}(\vec{e}))) \cong \mathbf{A}/\mathsf{Cg}^{\mathbf{A}}(\vec{0},\vec{e}) \\ \mathbf{P}_{1}\cong \mathbf{A}/\mathsf{Cg}^{\mathbf{A}}(\vec{1},\varphi_{\mathbf{A}}(\mu_{\mathbf{A}}(\vec{e})))\cong \mathbf{A}/\mathsf{Cg}^{\mathbf{A}}(\vec{1},\vec{e}).
\end{array}
\]
Therefore, for general reasons we get:

\[
\begin{array}{ccc}
\mathsf{Cg}^{\mathbf{A}}(\vec{0},\varphi_{\mathbf{A}}(\mu_{\mathbf{A}}(\vec{e})))=\mathsf{Cg}^{\mathbf{A}}(\vec{0},\vec{e}) \\ \mathsf{Cg}^{\mathbf{A}}(\vec{1},\varphi_{\mathbf{A}}(\mu_{\mathbf{A}}(\vec{e})))=\mathsf{Cg}^{\mathbf{A}}(\vec{1},\vec{e}).
\end{array}
\]
Hence, from Corollary 4 of \cite{V1999} it follows that $\varphi_{\mathbf{A}}(\mu_{\mathbf{A}}(\vec{e}))=\vec{e}$. On the other hand, if $g\in \V(\mathbf{0}\times \mathbf{0}, \mathbf{A})$ then we get $\mu_{\mathbf{A}}(\varphi_{\mathbf{A}}(g))=\text{!`}_{\vec{0},\varphi_{\mathbf{A}}(g)}\times \text{!`}_{\vec{1},\varphi_{\mathbf{A}}(g)}$. Thus we have $\mathbf{A}\cong \mathbf{A}/\mathsf{Cg}^{\mathbf{A}}(\vec{0},\varphi_{\mathbf{A}}(g))\times \mathbf{A}/\mathsf{Cg}^{\mathbf{A}}(\vec{1},\varphi_{\mathbf{A}}(g))$. Since $\V$ is coextensive, there exist unique $u:\mathbf{0}\rightarrow \mathbf{A}/\mathsf{Cg}^{\mathbf{A}}(\vec{0},\varphi_{\mathbf{A}}(g))$ and $v:\mathbf{0}\rightarrow \mathbf{A}/\mathsf{Cg}^{\mathbf{A}}(\vec{1},\varphi_{\mathbf{A}}(g))$ such that $g=u\times v$. Observe that due to $\mathbf{0}$ is initial in $\V$, it must be the case that $u=\text{!`}_{\vec{0},\varphi_{\mathbf{A}}(g)}$ and $v=\text{!`}_{\vec{1},\varphi_{\mathbf{A}}(g)}$ so $g=\mu_{\mathbf{A}}(\varphi_{\mathbf{A}}(g))$, as desired.  
\end{proof}
At this stage, one may be wondering if a characterization of coextensive varieties in terms of the representability of the functor $Z$ can be established. In \cite{Z2021} it was shown that not every variety with BFC and $\vec{0}$ and $\vec{1}$ is center stable and even if it is, it may be the case that it could be not coextensive. We conclude this section providing an effective answer to this question.

\begin{theorem}\label{converse coextensive}
Let $\V$ be a center stable variety. Then, the following are equivalent:
\begin{itemize}
\item[(1)] $\V$ is coextensive.
\item[(2)] The following conditions hold: 
\begin{itemize}
\item[(a)] The functor $Z:\V\rightarrow \Set$ is representable.
\item[(b)] The functor $\times: \mathbf{0}/\V\times \mathbf{0}/\V \rightarrow (\mathbf{0}\times \mathbf{0})/\V$ is full and faithful.
\end{itemize}
\end{itemize}
\end{theorem}
\begin{proof}
On the one hand, if we assume $(1)$, then $(2)(a)$ is in fact Theorem \ref{centrals bijection} and $(2)(b)$ follows from the definition of coextensivity. On the other hand, we start by showing that $(2)(a)$ implies that the functor $\times$ is essentially surjective in objects. To accomplish this, let $\mathbf{A} \in \mathcal{V}$, and consider $f: \mathbf{0} \times \mathbf{0} \rightarrow \mathbf{A}$. By assumption, there exists a bijection between $Z(\mathbf{A})$ and $\mathcal{V}(\mathbf{0} \times \mathbf{0}, \mathbf{A})$, allowing us to select $\vec{e}_{f} \in Z(\mathbf{A})$ corresponding to $f$ via this bijection. Take $\mathbf{A}_{0} = \mathbf{A}/\mathsf{Cg}^{\mathbf{A}}(\vec{e}_{f},\vec{0})$ and $\mathbf{A}_{1} = \mathbf{A}/\mathsf{Cg}^{\mathbf{A}}(\vec{e}_{f},\vec{1})$. Due to the coextensivity of $\mathcal{V}$, it is evident that $f$ is isomorphic to $\times (\text{!`}_{\vec{0},\vec{e}_{f}},\text{!`}_{\vec{1},\vec{e}_{f}})$, as claimed. Furthermore, by $(2)(b)$, $\times$ is an equivalence. Thus, $(1)$ holds, as desired.
\end{proof}

\section{The Gaeta topos and coextensive varieties}\label{The Gaeta topos and fp-coextensive varieties}

In this section we show that from the characterization of coextensive varieties described in Theorem \ref{charcoextensivity}, it is possible to provide a suitable axiomatization of the theory of indecomposable objects in the variety. Afterwards, we lift such a theory into the topos setting for defining the theory of central-free models. Subsequently, we show that the full subcategory of finitely presented algebras is closed under binary products. Finally, under the assumption that $\mathcal{V}$ is $(\vec{0},\vec{1})$-dense coextensive, we show that the Gaeta topos classifies central-free $\mathcal{V}$-models. \\

We start by proving some technical results on coextensive varieties which will be used along this section.  

Let $\mathcal{V}$ be a coextensive variety. Recall that again from Theorem \ref{charcoextensivity}, the relation $\vec{e}\diamond_{\mathbf{A}}\vec{f}$ in $\V$ is equationally definable. So we can assume that  
\[\sigma(\vec{x},\vec{y}):=\bigwedge_{i=1}^{k}p_{i}(\vec{x},\vec{y})=q_{i}(\vec{x},\vec{y})\]
defines the relation $\vec{e}\diamond_{\mathbf{A}}\vec{f}$ in $\V$. Now for the rest of this section we consider:
\begin{equation}\label{definincion teta}
\theta=\bigvee_{l=1}^{k}\mathsf{Cg}^{\mathbf{F}_{\V}(\vec{x},\vec{y})}(p_{l}(\vec{x},\vec{y}),q_{l}(\vec{x},\vec{y})),
\end{equation}
\[\mu=\mathsf{Cg}^{\mathbf{F}_{\V}(\vec{x},\vec{y})}(\vec{x},\vec{0})\vee \mathsf{Cg}^{\mathbf{F}_{\V}(\vec{x},\vec{y})}(\vec{y},\vec{1}),\] 
\[\lambda=\mathsf{Cg}^{\mathbf{F}_{\V}(\vec{x},\vec{y})}(\vec{x},\vec{1})\vee \mathsf{Cg}^{\mathbf{F}_{\V}(\vec{x},\vec{y})}(\vec{y},\vec{0}).\]

\begin{lemma}\label{pre useful lemma}
Let $\V$ be a coextensive variety. Then, the following hold:
\begin{itemize}
\item[(1)] $\theta\subseteq \mu, \lambda$;
\item[(2)] $\mu \circ \lambda =\lambda \circ \mu =\nabla^{\mathbf{F}_{\V}(\vec{x},\vec{y})}$;
\item[(3)]$\mathbf{F}_{\V}(\vec{x},\vec{y})/\mu \cong \mathbf{F}_{\V}(\vec{x},\vec{y})/\lambda$.
\end{itemize}
\end{lemma}
\begin{proof}

$(1)$ As we already know, $\V$ is a Pierce variety in which the relation $\vec{e}\diamond_{\mathbf{A}}{\vec{f}}$ is equationally definable, and in particular, it is a variety with $\vec{0}$ and $\vec{1}$.  Since $\vec{0}^{\mathbf{A}},\vec{1}^{\mathbf{A}}\in Z(\mathbf{A})$ for every $\mathbf{A}\in \V$, it is the case that
\[\V \models (\vec{x}=\vec{0} \wedge \vec{y}=\vec{1}) \rightarrow \sigma(\vec{x},\vec{y})\] 
so in particular
\[\V \models (\vec{x}=\vec{0} \wedge \vec{y}=\vec{1}) \rightarrow p_{l}(\vec{x},\vec{y})=q_{l}(\vec{x},\vec{y})\]
for every $1\leq l\leq k$. Then from Lemma \ref{vey useful lemma} 
\[\mathsf{Cg}^{\mathbf{F}_{\V}(\vec{x},\vec{y})}(p_{l}(\vec{x},\vec{y}),q_{l}(\vec{x},\vec{y}))\subseteq \mu\]
for every $1\leq l\leq k$. Therefore, $\theta \subseteq \mu$ as required. The proof of $\theta \subseteq \lambda$ is analogue.

$(2)$ Let $(s(\vec{x},\vec{y}),t(\vec{x},\vec{y}))\in \nabla^{\mathbf{F}_{\V}(\vec{x},\vec{y})}$. By Theorem \ref{charcoextensivity} (2) $\mathcal{V}$ has a decomposition term $U(x,y,\vec{z},\vec{w})$. Let us consider $r(\vec{x},\vec{y}):=U(s(\vec{x},\vec{y}),t(\vec{x},\vec{y}),\vec{x},\vec{y})$. Observe that $r(\vec{0},\vec{1})=s(\vec{0},\vec{1})$ and $r(\vec{1},\vec{0})=t(\vec{1},\vec{0})$. Thus, it follows that 
\[\V\models (\vec{x}=\vec{0} \wedge \vec{y}=\vec{1}) \rightarrow r(\vec{x},\vec{y})=s(\vec{x},\vec{y}) \]
and
\[\V\models (\vec{x}=\vec{1} \wedge \vec{y}=\vec{0}) \rightarrow r(\vec{x},\vec{y})=t(\vec{x},\vec{y}). \]
Thus by Lemma \ref{vey useful lemma} we get
\[\begin{array}{ccc}
(s(\vec{x},\vec{y}),r(\vec{x},\vec{y}))\in \mu & \text{and} & (r(\vec{x},\vec{y}),t(\vec{x},\vec{y}))\in \lambda.
\end{array}\]
Hence $\mu\circ \lambda = \nabla^{\mathbf{F}_{\V}(\vec{x},\vec{y})}$. Finally, we stress that if we take \[r'(\vec{x},\vec{y})=U(s(\vec{x},\vec{y}),t(\vec{x},\vec{y}),\vec{y},\vec{x})\] it is no hard to see that $r'(\vec{0},\vec{1})= t(\vec{0},\vec{1})$ and $r'(\vec{1},\vec{0})= s(\vec{1},\vec{0})$. Therefore, by the same argument we employed before, together with Lemma \ref{vey useful lemma}, we get $\lambda\circ \mu = \nabla^{\mathbf{F}_{\V}(\vec{x},\vec{y})}$, as claimed. 
\\
$(3)$ This follows from the fact that the map $g:\mathbf{F}_{\V}(\vec{x},\vec{y})/\mu \to \mathbf{F}_{\V}(\vec{x},\vec{y})/\lambda$, defined by $g(s(\vec{x},\vec{y})/\mu)=s(\vec{y},\vec{x})/\lambda$ is an isomorphism. The details are left to the reader.
\end{proof}

\begin{theorem}\label{decomposition theorem}
    In every coextensive variety $\mathcal{V}$:
    \[\mathbf{F}_{\V}(\vec{x},\vec{y})/\theta\cong \mathbf{F}_{\V}(\vec{x},\vec{y})/\mu \times \mathbf{F}_{\V}(\vec{x},\vec{y})/\lambda.\]
\end{theorem}
\begin{proof}
    Observe that according Lemmas \ref{factor congruence quotinents} and \ref{pre useful lemma}, it is enough to show that $\theta=\mu \cap \lambda$. For the sake of readability, along this proof, we will write $\mathbf{B}:=\mathbf{F}_{\V}(\vec{x},\vec{y})/\theta$, $\alpha:= \mathsf{Cg}^{\mathbf{B}}(\vec{x}/\theta,\vec{0}/\theta)$ and $\beta:= \mathsf{Cg}^{\mathbf{B}}(\vec{y}/\theta,\vec{0}/\theta)$. We start by noticing that $\vec{x}/\theta\diamond_{\mathbf{B}}\vec{y}/\theta$. Then, $\mathbf{B}\cong \mathbf{B}/\alpha \times \mathbf{B}/\beta$. Now, we consider the following morphisms:
    \begin{displaymath}
    \begin{array}{cc}
        \xymatrix{
        \mathbf{F}_{\V}(\vec{x},\vec{y}) \ar[r]^-{f_{\theta}} & \mathbf{B} \ar[r]^-{g_{\alpha}} & \mathbf{B}/\alpha
        } & \xymatrix{
        \mathbf{F}_{\V}(\vec{x},\vec{y}) \ar[r]^-{f_{\theta}} & \mathbf{B} \ar[r]^-{g_{\beta}} & \mathbf{B}/\beta
        }
    \end{array}
    \end{displaymath}
    Let $\gamma=\mathsf{Ker}(g_{\alpha}f_{\theta})$ and $\delta=\mathsf{Ker}(g_{\beta}f_{\theta})$, 
    By the Correspondence Theorem (Theorem 3.20 \cite{BS1981}), it is the case that $\alpha=\gamma/\theta$ and $\beta=\delta/\theta$.We proceed to show that $\gamma=\mu$. The proof that $\delta=\lambda$ is similar. To do so, observe that since $x_l/\theta=0_l/\theta$, for every $1\leq l\leq N$, it is readily seen that $g_{\alpha}f(x_l)=g_{\alpha}f(0_l)$. Consequently, $\mathsf{Cg}^{\mathbf{F}_{\V}(\vec{x},\vec{y})}(\vec{x},\vec{0})\subseteq \gamma$. Furthermore, since $\vec{y}/\theta$ is the complement of $\vec{x}/\theta$ in $Z(\mathbf{B})$, we get that $\alpha=\mathsf{Cg}^{\mathbf{B}}(\vec{y}/\theta,\vec{1}/\theta)$. From the latter, it is straightforward to see that $g_{\alpha}f(y_l)=g_{\alpha}f(1_l)$, for every $1\leq l\leq N$. Hence $\mathsf{Cg}^{\mathbf{F}_{\V}(\vec{x},\vec{y})}(\vec{y},\vec{1})\subseteq \gamma$. Now, let us assume that there exists $\varepsilon\in \mathsf{Con}(\mathbf{F}_{\V}(\vec{x},\vec{y}))$ such that $\mathsf{Cg}^{\mathbf{F}_{\V}(\vec{x},\vec{y})}(\vec{x},\vec{0}), \mathsf{Cg}^{\mathbf{F}_{\V}(\vec{x},\vec{y})}(\vec{y},\vec{1})\subseteq \varepsilon$. We need to see that $\gamma \subseteq \varepsilon$. Suppose the contrary. Then, there are $s(\vec{x},\vec{y}), t(\vec{x},\vec{y})\in \mathbf{F}_{\V}(\vec{x},\vec{y})$ such that (a) $(s(\vec{x},\vec{y}), t(\vec{x},\vec{y}))\in \gamma$ but (b) $(s(\vec{x},\vec{y}), t(\vec{x},\vec{y}))\notin \varepsilon$. Let $\mathbf{C}:=\mathbf{B}/\alpha \in \mathcal{V}$. Then, from (b), the assumption on $\varepsilon$ and Lemma \ref{vey useful lemma}, we have that the following hold:
    \begin{gather*}
        \mathcal{V}\models \vec{x}=\vec{0}\; \wedge\; s(\vec{x}, \vec{y})\neq t(\vec{x}, \vec{y})\\
        \mathcal{V}\models \vec{y}=\vec{1}\; \wedge\; s(\vec{x}, \vec{y})\neq t(\vec{x}, \vec{y})
    \end{gather*}
    Hence, 
    \begin{gather*}
        \mathcal{V}\nvDash (\vec{x}=\vec{0}\; \wedge\; \vec{y}=\vec{1})\; \to  s(\vec{x}, \vec{y})= t(\vec{x}, \vec{y}).
    \end{gather*}
    But from (a) and the fact that $\vec{x}/\theta=\vec{0}/\theta$ and $\vec{y}/\theta=\vec{1}/\theta$ it follows that $s^{\mathbf{C}}(\vec{0}^{\mathbf{C}}, \vec{1}^{\mathbf{C}})=t^{\mathbf{C}}(\vec{0}^{\mathbf{C}}, \vec{1}^{\mathbf{C}})$, which is absurd. Then it must be that $\gamma\subseteq \varepsilon$. Therefore, $\gamma=\mu$, as desired. Thus, by Lemma \ref{factor congruence quotinents}, $\gamma \circ \delta =\delta \circ \gamma =\nabla^{\mathbf{F}_{\V}(\vec{x},\vec{y})}$  and $\theta=\gamma \cap \delta$. So, based on what we just have proved, we conclude that $\theta=\mu \cap \lambda$. Finally, from lemmas \ref{pre useful lemma} and \ref{factor congruence quotinents}, we are able to conclude that $\mathbf{F}_{\V}(\vec{x},\vec{y})/\theta\cong \mathbf{F}_{\V}(\vec{x},\vec{y})/\mu\times \mathbf{F}_{\V}(\vec{x},\vec{y})/\lambda$, as claimed.
\end{proof}

Now we will prove a much more general result concerning coextensive varieties. It states that binary products of free finitely generated algebras is finitely presented. We start by establishing some notation. Let $n,m\in \mathbb{N}$, fro the following we set $n+m:=\{x_1,\ldots, x_n,y_1\ldots, y_m\}$.

\begin{lemma}\label{prodfg is fingen}
Let $\mathbf{P}$ be the product of $\mathbf{F}_{\V}(n)$ and $\mathbf{F}_{\V}(m)$. Then $\mathbf{P}$ is finitely generated.
\end{lemma}
\begin{proof}
    Let $c$ be a fixed constant symbol in the language of $\mathcal{V}$ and let $r\in \mathbb{N}_{>0}$. We will write $\vec{c}_{r}$ for the $r$-tuple whose coordinates are always $c$. Furthermore, we set $\vec{x}:=(x_1,\ldots,x_n)$ and $\vec{y}:=(y_1,\ldots,y_m)$. We claim that $\mathbf{P}$ is generated by the set $S=\{[\vec{x},\vec{c}_n], [\vec{c}_m,\vec{y}],[\vec{0},\vec{1}], [\vec{1},\vec{0}]\}$. Indeed, let $(r(\vec{x}),v(\vec{y}))\in P$ and consider 
    \[\psi(\vec{x},\vec{y},\vec{z},\vec{w})=U(r(\vec{x}), v(\vec{y}), \vec{z},\vec{w}),\]
    where $U$ is the decomposition term of $\mathcal{V}$. Observe that from the definition of $U$, it is readily seen that \[\psi^{\mathbf{P}}([\vec{x},\vec{c}_n], [\vec{c}_m,\vec{y}],[\vec{0},\vec{1}], [\vec{1},\vec{0}])=(r(\vec{x}),v(\vec{y})).\]
\end{proof}

Let $\vec{z}:=(z_1,\ldots, z_N)$ and $\vec{w}:=(w_1,\ldots, w_N)$. We write $n+m+\{\vec{z}, \vec{w}\}$ to denote the set $\{x_1,\ldots, x_n,y_1\ldots, y_m, z_1,\ldots, z_N, w_1,\ldots, w_N\}$. Consider now the homomorphisms $L, M:\mathbf{F}_{\V}(n+m+\{\vec{z}, \vec{w}\})\to \mathbf{P}$ which are determined by assignments as described below:
\begin{displaymath}
    \begin{array}{ccc}
         &  & \\
        L(x_i)=x_i, &  & M(x_i)=c, \\
        L(y_j)=c, &  & M(y_j)=y_i, \\
        L(z_l)=0_l, &  & M(z_l)=1_l, \\
        L(w_l)=1_l, &  & M(w_l)=0_l, \\
    \end{array}
\end{displaymath}
where $1\leq i\leq n$, $1\leq j\leq m$ and $1\leq l\leq N$.
\begin{lemma}\label{map surjective}
    The map \[\langle L,M\rangle:\mathbf{F}_{\V}(n+m+\{\vec{z}, \vec{w}\})\to \mathbf{P}\] is surjective.
\end{lemma}
\begin{proof}
    Immediate from Lemma \ref{prodfg is fingen}.
\end{proof}

\begin{lemma}\label{generators kernel}
    In $\mathbf{F}_{\V}(n+m+\{\vec{z}, \vec{w}\})$, the following hold:
    \begin{enumerate}
        \item For every $1\leq j\leq m$, \[\langle L,M\rangle(U(y_j,c,\vec{z},\vec{1}))=\langle L,M\rangle(U(y_j,c,\vec{1},\vec{z})).\] 
        \item For every $1\leq i\leq n$, \[\langle L,M\rangle(U(x_i,c,\vec{w},\vec{1}))=\langle L,M\rangle(U(x_i,c,\vec{1},\vec{w})).\]
        \item For every $1\leq l\leq k$,
        \[\langle L,M\rangle(p_l(\vec{z},\vec{w}))=\langle L,M\rangle(q_l(\vec{z},\vec{w})).\]
    \end{enumerate}
\end{lemma}
\begin{proof}
    (1) Let $1\leq j\leq m$. Then we have 
    \begin{displaymath}
        \begin{array}{rcl}
            \langle L,M\rangle(U(y_j,c,\vec{z},\vec{1})) & = & U^{\mathbf{P}}((c,y_j),(c,c),[\vec{0},\vec{1}], [\vec{1},\vec{1}]) \\
             & = & (U(c,c,\vec{0},\vec{1}),U(y_j,c,\vec{1},\vec{1})) \\
             & = & (c,U(y_j,c,\vec{1},\vec{1})) \\
             & = & (U(c,c,\vec{1},\vec{0}),U(y_j,c,\vec{1},\vec{1}))\\
             & = & \langle L,M\rangle(U(y_j,c,\vec{1},\vec{z})).
        \end{array}
    \end{displaymath}
    (2) Analogue to (1).
    \\
    (3) Notice that for every $1\leq l\leq k$:
    \begin{displaymath}
        \begin{array}{c}
           \langle L,M\rangle (p_l(\vec{z},\vec{w}))=(p^{{F}_{\V}(n)}_{l}(\vec{0},\vec{1}), 
p^{{F}_{\V}(m)}_{l}(\vec{1},\vec{0})), \\
           \langle L,M\rangle (q_l(\vec{z},\vec{w}))=(q^{{F}_{\V}(n)}_{l}(\vec{0},\vec{1}), 
q^{{F}_{\V}(m)}_{l}(\vec{1},\vec{0})). 
        \end{array}
    \end{displaymath}   
Since $\vec{0}\diamond_{\mathbf{F}_{\V}(n)}\vec{1}$ and  $\vec{1}\diamond_{\mathbf{F}_{\V}(m)}\vec{0}$, for every $1\leq l\leq k$ we get:
    \begin{displaymath}
        \begin{array}{ccc}
            p^{{F}_{\V}(n)}_{l}(\vec{0},\vec{1})= q^{{F}_{\V}(n)}_{l}(\vec{0},\vec{1}),& & p^{{F}_{\V}(m)}_{l}(\vec{1},\vec{0})= q^{{F}_{\V}(m)}_{l}(\vec{1},\vec{0}).
        \end{array}
    \end{displaymath}
    This allows us to conclude that $\langle L,M\rangle (p_l(\vec{z},\vec{w}))=\langle L,M\rangle (q_l(\vec{z},\vec{w}))$, as desired.
\end{proof}

\begin{theorem}\label{fine theorem}
Binary products of free finitely generated algebras is finitely presented.    
\end{theorem}
\begin{proof}
Let $\mathbf{D}=\mathbf{F}_{\V}(n+m+\{\vec{z}, \vec{w}\})$ and let us consider the following sets of $D^2$:  
\[A=\{(U(y_j,c,\vec{z},\vec{1}),U(y_j,c,\vec{1},\vec{z}))\colon 1\leq j\leq m \},\]
\[B=\{(U(x_i,c,\vec{w},\vec{1}),U(x_i,c,\vec{1},\vec{w}))\colon 1\leq i\leq m \},\]
\[C=\{(p_l(\vec{z}, \vec{w}),q_l(\vec{z}, \vec{w}))\colon 1\leq l\leq k\}.\]

Let $\chi$ be the least congruence of $\mathbf{D}$ generated by $A\cup B\cup C$. We will show that $\mathbf{D}/\chi$ is isomorphic to $\mathbf{P}$. From Lemma \ref{generators kernel}, there exists a unique homomorphism $\Gamma:\mathbf{D}/\chi\to \mathbf{P}$ such that the following diagram 
\begin{displaymath}
    \xymatrix{
    \mathbf{D}\ar[r] \ar[dr]_-{\langle L, M\rangle} & \mathbf{D}/\chi \ar@{-->}[d]^-{\Gamma}\\
     & \mathbf{P}
    }
\end{displaymath}
commutes. Moreover, from Lemma \ref{map surjective}, $\Gamma$ is surjective. It only remains to show that $\Gamma$ is injective. To do so, let  $t(\vec{x},\vec{y}, \vec{z},\vec{w}),s(\vec{x},\vec{y}, \vec{z},\vec{w})\in D$ such that 
\[\Gamma(
t(\vec{x},\vec{y}, \vec{z},\vec{w})/\chi)=\Gamma(
s(\vec{x},\vec{y}, \vec{z},\vec{w})/\chi).\] 
Then we get 
\[t^{\mathbf{P}}((\vec{x},\vec{c}_n),(\vec{c}_m,\vec{y}),[\vec{0},\vec{1}], [\vec{1},\vec{0}])=s^{\mathbf{P}}((\vec{x},\vec{c}_n),(\vec{c}_m,\vec{y}),[\vec{0},\vec{1}], [\vec{1},\vec{0}]),\]
hence:
\begin{equation}\label{key equation1}
t(\vec{x},\vec{c}_m,\vec{0},\vec{1})=s(\vec{x},\vec{c}_m,\vec{0},\vec{1}),
\end{equation}
\begin{equation}\label{key equation2}
t(\vec{c}_n,\vec{y},\vec{1},\vec{0})=s(\vec{c}_n,\vec{y},\vec{1},\vec{0}).
\end{equation}
In order to see that $(t(\vec{x},\vec{y}, \vec{z},\vec{w}),s(\vec{x},\vec{y}, \vec{z},\vec{w}))\in \chi$, we will use Lemma \ref{vey useful lemma}. Therefore, we need to show that 
\[\mathcal{V}\models \phi_1(\vec{y},\vec{z})\wedge \phi_2(\vec{x},\vec{w}) \wedge \phi_1(\vec{z},\vec{w})\to t(\vec{x},\vec{y}, \vec{z},\vec{w})=s(\vec{x},\vec{y}, \vec{z},\vec{w}),\]
where
\begin{equation}\label{ecua1} \phi_1(\vec{y},\vec{z}):=\bigwedge_{j=1}^{m}U(y_j,c,\vec{z},\vec{1})=U(y_j,c,\vec{1},\vec{z}),   
\end{equation}
\begin{equation}\label{ecua2}   \phi_2(\vec{x},\vec{w}):=\bigwedge_{i=1}^{m}U(x_i,c,\vec{w},\vec{1})=U(x_i,c,\vec{1},\vec{w}),
\end{equation}
\begin{equation}\label{ecua3}
\phi_3(\vec{z},\vec{w}):=\bigwedge_{l=1}^{k}p_{l}(\vec{z},\vec{w})=q_{l}(\vec{z},\vec{w}).
\end{equation}
Let $\mathbf{A}\in \mathcal{V}$ and let $\alpha:n+m+\{\vec{z}, \vec{w}\}\to A$ be an assignment. By (\ref{ecua3}) it is the case that $\alpha(\vec{z})\diamond_{\mathbf{A}}\alpha(\vec{w})$, so $\mathbf{A}\cong \mathbf{A}/\theta_{\vec{z}}\times \mathbf{A}/\theta_{\vec{w}}$, where $\theta_{\vec{z}}=\mathsf{Cg}^{\mathbf{A}}(\alpha(\vec{z}),\vec{0})$ and $\theta_{\vec{w}}=\mathsf{Cg}^{\mathbf{A}}(\alpha(\vec{w}),\vec{0})$. Therefore, $t^{\mathbf{A}}(\alpha(\vec{x}),\alpha(\vec{y}),\alpha(\vec{z}),\alpha(\vec{w}))$ can be expressed as the ordered pair whose first coordinate is 
\begin{equation}\label{term1}
t^{\mathbf{A}/\theta_{\vec{z}}}(\alpha(\vec{x})/\theta_{\vec{z}},\alpha(\vec{y})/\theta_{\vec{z}},\alpha(\vec{z})/\theta_{\vec{z}},\alpha(\vec{w})/\theta_{\vec{z}})
\end{equation}
and its second coordinate is 
\begin{equation}\label{term2}
t^{\mathbf{A}/\theta_{\vec{w}}}(\alpha(\vec{x})/\theta_{\vec{w}},\alpha(\vec{y})/\theta_{\vec{w}},\alpha(\vec{z})/\theta_{\vec{w}},\alpha(\vec{w})/\theta_{\vec{w}}).
\end{equation}
Observe now, that from (\ref{ecua1}), (\ref{ecua2}) and (\ref{Ecucentrals}), we have that $\alpha(\vec{y})/\theta_{\vec{z}}=\vec{c}_{m}^{\mathsf{A}}/\theta_{\vec{z}}$ and $\alpha(\vec{x})/\theta_{\vec{w}}=\vec{c}_{n}^{\mathsf{A}}/\theta_{\vec{w}}$, and since $\alpha(\vec{z})/\theta_{\vec{z}}=\vec{0}/\theta_{\vec{z}}$, $\alpha(\vec{w})/\theta_{\vec{z}}=\vec{1}/\theta_{\vec{z}}$, $\alpha(\vec{z})/\theta_{\vec{w}}=\vec{1}/\theta_{\vec{w}}$ and $\alpha(\vec{w})/\theta_{\vec{w}}=\vec{0}/\theta_{\vec{w}}$, then (\ref{term1}) and (\ref{term2}) take the form 
\[
\alpha(t(\vec{x},\vec{c}_{m},\vec{0},\vec{1}))/\theta_{\vec{z}}
\]
and
\[
\alpha(t(\vec{c}_{n},\vec{y},\vec{1},\vec{0}))/\theta_{\vec{w}},
\]
respectively. From the latter, (\ref{key equation1}) and (\ref{key equation2}), it is no hard to see that one obtains that $$t^{\mathbf{A}}(\alpha(\vec{x}),\alpha(\vec{y}),\alpha(\vec{z}),\alpha(\vec{w}))=s^{\mathbf{A}}(\alpha(\vec{x}),\alpha(\vec{y}),\alpha(\vec{z}),\alpha(\vec{w})).$$ 
Hence, by Lemma \ref{vey useful lemma} we get that $t(\vec{x},\vec{y}, \vec{z},\vec{w})/\chi=s(\vec{x},\vec{y}, \vec{z},\vec{w})/\chi$, as expected. So  $\Gamma$ is injective and consequently, an isomorphism. Finally, since $\chi$ is finitely generated in $\mathbf{D}$, we are able to coclude that $\mathbf{P}$ is finitely presented, as desired. 
\end{proof}

We emphasize that, although Theorem \ref{decomposition theorem} allows us to decompose the quotient $\mathbf{F}_{\V}(\vec{x},\vec{y})/\theta$ in terms of the congruences $\lambda$ and $\mu$, our main interest lies in those coextensive varieties where each factor is isomorphic to $\mathbf{0}$. This will be a key point in the following section. While this holds in many known cases —such as distributive lattices and commutative rings with unit— it does not hold in general. To illustrate this, let $\mathcal{R}$ be the coextensive variety of commutative rings with unit, and consider the class $\mathcal{K} = \{(\mathbf{A}, c^{\mathbf{A}}) : \mathbf{A} \in \mathcal{R}\}$, where $c$ is a new constant symbol distinct from $0$ and $1$. It is not difficult to see that $\mathcal{K}$ is a coextensive variety; however, 
\[\mathbf{F}_{\mathcal{K}}(\emptyset) = \{0, 1\} \cup \{p(\vec{c}) : p(\vec{x}) \in T_{\mathcal{R}}(n),; n \in \mathbb{N}\}\neq \mathsf{Sg}^{\mathbf{F}_{\mathcal{K}}(\emptyset)}(0,1).\]
Therefore, in general, $\mathbf{F}_{\mathcal{V}}(\emptyset)$ is not generated by ${\vec{0}, \vec{1}}$. This observation motivates the following definition: 
\begin{definition}\label{01-dense}
    A coextensive variety $\mathcal{V}$ is said to be $(\vec{0},\vec{1})$-dense if $\mathbf{0}$ is generated by $\vec{0}$ and $\vec{1}$. 
\end{definition} 
We conclude this section by presenting a result that demonstrates that Definition \ref{01-dense} captures the precise notion we need for our purposes:
\begin{lemma}\label{char 01-dense}
  Let $\mathcal{V}$ be a coextensive variety. Then, the following are equivalent:
  \begin{enumerate}
      \item $\mathcal{V}$ is $(\vec{0},\vec{1})$-dense.
      \item For every $n$-ary term $p(\vec{z})$ and constant symbols $c_{1},...,c_{n}$ in the language of $\V$, there exists a $2N$-ary term $q(\vec{x},\vec{y})$ such that 
\[\V \models p(\vec{c}) = q(\vec{0},\vec{1})\]
where $\vec{c}$ stands for the n-tuple $(c_1,\ldots, c_n)$.
\item $\mathbf{F}_{\V}(\vec{x},\vec{y})/\mu \cong \mathbf{0} \cong \mathbf{F}_{\V}(\vec{x},\vec{y})/\lambda $.
  \end{enumerate}
\end{lemma}
\begin{proof}
    (1)$\Rightarrow$(2): Let $\mathbf{A}\in \mathcal{V}$, let $p(\vec{z})$ be an $n$-ary term and let $c_{1},...,c_{n}$ be constant symbols in the language of $\V$. By (1), there exists $q(\vec{x}, \vec{y})$ such that $p(\vec{c})=
q(\vec{0},\vec{1})$. Since $\mathbf{F}_{\V}(\emptyset)\cong \mathbf{0}$, then, it is the case that \[p^{\mathbf{A}}(\vec{c}^{\mathbf{A}})=\text{!`}_{\mathbf{A}}(p(\vec{c}))=\text{!`}_{\mathbf{A}}(q(\vec{0},\vec{1}))= q^{\mathbf{A}}(\vec{0}^{\mathbf{A}},\vec{1}^{\mathbf{A}}).\]
So (2) holds.
\\
(2)$\Rightarrow$(3): Let $f:\mathbf{F}_{\V}(\vec{x},\vec{y})\rightarrow \mathbf{F}_{\V}(\vec{x},\vec{y})/\mu$ be the unique homomorphism determined by the assignments $f(\vec{x})=\vec{0}$ and $f(\vec{y})=\vec{1}$. Observe that from (2), $f$ is surjective. Since $(\vec{x}, \vec{0}), (\vec{0}, \vec{1})\in \mathsf{Ker}(f)$, there exists a unique $\Gamma$ making the following diagram
\begin{displaymath}
    \xymatrix{
    \mathbf{F}_{\V}(\vec{x},\vec{y})\ar[r] \ar[dr]_-{\langle L, M\rangle} & \mathbf{F}_{\V}(\vec{x},\vec{y})/\mu \ar@{-->}[d]^-{\Gamma}\\
     & \mathbf{P}
    }
\end{displaymath}
commutes. Moreover, $\Gamma$ is surjective as well, since $f$ is. It only remains to show the injectivity of $\Gamma$. So supposse that $\Gamma(s(\vec{x}, \vec{y})/\mu)=\Gamma(t(\vec{x}, \vec{y})/\mu)$. Then $s(\vec{0}, \vec{1})=t(\vec{0}, \vec{1})$. From this, it follows easily that $\mathcal{V}\models \vec{x}=\vec{0} \wedge \vec{y}=1 \rightarrow s(\vec{x}, \vec{y})=t(\vec{x}, \vec{y})$. Hence, by Lemma \ref{vey useful lemma}, $s(\vec{x}, \vec{y})/\mu=t(\vec{x}, \vec{y})/\mu$. This shows that $\Gamma$ is an isomorphism. The proof that $\mathbf{F}_{\V}(\vec{x},\vec{y})/\lambda$ is isomorphic to $\mathbf{0}$ is analogue.
\\
(3)$\Rightarrow$(1): Let $p(\vec{z})$ be an $n$-ary term and let $c_{1},...,c_{n}$ be constant symbols in the language of $\V$. By (3) there is an isomorphism $\tau:\mathbf{F}_{\V}(\vec{x},\vec{y})/\mu \to \mathbf{0}$. If we consider the composite $\tau q_{\mu}:\mathbf{F}_{\V}(\vec{x},\vec{y}) \to \mathbf{0}$, which is clearly surjective, then we obtain: 
\begin{displaymath}
    \begin{array}{rcll}
        p(\vec{c}) & = &  \tau q_{\mu}(q(\vec{x},\vec{y})), & \text{for some}\; q(\vec{x},\vec{y})\in F_{\V}(\vec{x},\vec{y}). \\
         & = & q(\tau q_{\mu}(\vec{0}),\tau q_{\mu}(\vec{1}))) & \\
         & = & q(\vec{0},\vec{1}). &
    \end{array}
\end{displaymath}
Therefore, (1) holds. This concludes the proof.
\end{proof}
As an straightforward consequence of lemmas \ref{char 01-dense} and \ref{pre useful lemma} we get the following:
\begin{corollary}\label{decomposition of 0x0}
    In every $(\vec{0},\vec{1})$-dense coextensive variety $\mathcal{V}$:
    \[
    \mathbf{0}\times \mathbf{0} \cong \mathbf{F}_{\V}(\vec{x},\vec{y})/\mu \times \mathbf{F}_{\V}(\vec{x},\vec{y})/\lambda\cong \mathbf{F}_{\V}(\vec{x},\vec{y})/\theta.
    \]
\end{corollary}

\subsection{Central-free models in a topos}\label{Indecomposable models in a Topos}

Let $\V$ be a coextensive variety. An algebra $\mathbf{A}$ of $\mathcal{V}$ is said to be \emph{directly indecomposable} if is not trivial and, moreover, if $\mathbf{A}\cong \mathbf{B}\times \mathbf{C}$, then $\mathbf{B}$ or $\mathbf{C}$ is trivial.  The following result allows to show that from the theory of central elements it is possible to find an axiomatization for the theory of directly indecomposable objects of $\mathcal{V}$.

\begin{lemma}\label{Theory of connected models}
In every coextensive variety $\V$, the class of directly indecomposable objects is axiomatizable by a first order formula.
\end{lemma}
\begin{proof}
From Theorem \ref{charcoextensivity} (2) the relation $\vec{e}\diamond_{\mathbf{A}} \vec{f}$ is equationally definable in $\mathcal{V}$. So we can take $\sigma(\vec{x},\vec{y})$ as a finite conjunction of equations in the variables $\vec{x},\vec{y}$ defining such a relation. It is immediate that $\mathbf{A}\in \V$ is directly indecomposable if and only if in $\mathbf{A}$ the following sentence holds
\begin{center}
$\vec{0}\neq \vec{1}$ and $(\forall \vec{e},\vec{f}\; \sigma(\vec{e},\vec{f})\rightarrow ((\vec{e}=\vec{0} \wedge \vec{f}=\vec{1})\vee (\vec{e}=\vec{1} \wedge \vec{f}=\vec{0}))).$ 
\end{center} 
\end{proof}

Let $\V$ be a coextensive variety and let $\mathsf{E}$ be a topos. We recall that an object $M$ in $\mathsf{E}$ is a \emph{$\mathcal{V}$-model} if all the equations defining $\V$ (expressed as suitable sequents) hold in $M$. The $\mathcal{V}$-models in $\mathsf{E}$ with their morphisms (i.e. the arrows in $\mathsf{E}$ that preserve the structure of the models) form naturally a category $\mathsf{Mod}(\V,\mathsf{E})$. It is clear that, when regarding the topos $\Set$, the category $\mathsf{Mod}(\V,\mathsf{Set})$ coincides with $\V$ (for details, the reader may consult section D1 of \cite{J2002}).
\\

Now, inspired on Lemma \ref{Theory of connected models} we introduce the theory to whose models, we will devote our attention throughout Section \ref{The characterization}.

\begin{definition}\label{indecomposable sequents}
Let $\mathsf{E}$ be a topos. A $\mathcal{V}$-model $M$ of $\mathsf{E}$ is central-free if the sequents
\[0=1 \vdash \perp\]
\[\sigma(\vec{x},\vec{y})\vdash_{\vec{x},\vec{y}} (\vec{x}=\vec{0} \wedge \vec{y}=\vec{1}) \vee(\vec{x}=\vec{1} \wedge \vec{y}=\vec{0})\]
hold in the internal logic of $\mathsf{E}$.
\end{definition}

We stress that Definition \ref{indecomposable sequents} is addressed to axiomatize those $\mathcal{V}$-models in a topos which cannot be ``decomposed'' by the employment of complementary central elements. Observe that by Lemma \ref{Theory of connected models}, in the topos $\Set$, central-free models, coincide with directly indecomposable algebras of $\V$. Nevertheless, we are aware that, at this point, the reader may be wondering why we choose the name central-free rather than directly indecomposable for the theory just introduced. We must declare that this was on purpose, and is justified by the fact that although there is evidence of cases in which in a topos, directly indecomposable $\V$-models coincide with central-free $\V$-models (as the case of the theory of bounded distributive lattices which was proved in \cite{Za2017}), in our setting, such a fact is not relevant for the r\^{o}le that central-free models will take in the rest of this paper.
\\

Next, we show that Definition \ref{indecomposable sequents} can be easily internalized (in the sense of translating it in terms of certain diagrams in $\mathsf{E}$) with the application of very basic tools from the internal logic of toposes.

\begin{lemma}\label{Connected VModels in a topos}
Let $\mathsf{E}$ be a topos and let $M$ be a $\V$-model in $\mathsf{E}$. The following are equivalent:
\begin{itemize}
\item[(1)] $M$ is central-free,
\item[(2)] The diagram below 
\begin{displaymath}
\xymatrix@1{
0 \ar[r]^-{\text{!`}_{1}}  & 1 \ar@<0.5ex>[r]^-{\vec{1}_M} \ar@<-0.5ex>[r]_-{\vec{0}_M}  & M^{N} 
}
\end{displaymath}
is an equalizer in $\mathsf{E}$, and the morphism $\alpha: 1+1\rightarrow [\sigma(\vec{x},\vec{y})]_{M}$ is an isomorphism.
\end{itemize}
\end{lemma}
\begin{proof}
Let $M$ be a $\V$-model in $\mathsf{E}$ and let $\vec{0}_{M}$, $\vec{1}_{M}$ and $[\sigma(\vec{x},\vec{y})]_{M}$ be the interpretations in $M$ of the constants $\vec{0}$, $\vec{1}$ and of an equation $\sigma(\vec{x},\vec{y})$ defining the relation $\vec{e}\diamond_{\mathbf{A}}\vec{f}$ in $\V$, respectively. It is clear from the fact that the interpretation of $\perp$ is the initial object in $\mathsf{E}$, that the first diagram is an equalizer if and only if the first sequent in Definition \ref{indecomposable sequents} hold. On the other hand, notice that $1+1$ coincides with the interpretation in $M$ of the formula $(\vec{x}=\vec{0} \wedge \vec{y}=\vec{1}) \vee(\vec{x}=\vec{1} \wedge \vec{y}=\vec{0})$, so let us consider the elements $\langle \vec{0}_{M},\vec{1}_{M}\rangle: 1\rightarrow M^{N}\times M^{N}$ and $\langle \vec{1}_{M},\vec{0}_{M}\rangle: 1\rightarrow M^{N}\times M^{N}$. If $\alpha=[\langle \vec{0}_{M},\vec{1}_{M} \rangle, \langle \vec{1}_{M},\vec{0}_{M}\rangle]$ denotes the morphism from $1+1$ to $[\sigma(\vec{x},\vec{y})]_{M}$ induced by the coproduct, then it is no hard to see that $\alpha$ is an isomorphism if and only if the second sequent in Definition \ref{indecomposable sequents} hold. This concludes the proof.	 

\end{proof}

\subsection{The characterization}\label{The characterization}

In every small extensive category $\mathsf{C}$, finite families $\{\;f_{i}\colon X_{i}\rightarrow X\mid i\in I\;\}$ such that the coproduct induced arrow $\Sigma X_{i}\rightarrow X$ is an isomorphism form the basis of a Grothendieck topology. The topology $J_{\mathcal{G}}$ generated by such a basis is called the \emph{Gaeta Topology} and the \emph{Gaeta topos} $\mathcal{G}(\mathsf{C})$, is the topos of sheaves on the site $(\mathsf{C}, J_{\mathcal{G}})$. As observed in \cite{CPR2001},  $\mathcal{G}(\mathsf{C})$ is equivalent to the category $\mathsf{Lex}(\mathsf{C}^{\mathrm{op}},\Set)$ of product preserving functors to $\Set$ from the category $\mathsf{C}^{\mathrm{op}}$ with finite products. This fact implies that $J_{\mathcal{G}}$ is subcanonical (i.e.  all representable presheaves on this site are sheaves).

\begin{remark}\label{Gaeta continous}
Let $\mathsf{C}$ be an extensive category with a terminal object $1$ and let $\mathsf{E}$ be a topos. Notice that a finite limit preserving functor $G:\mathsf{C} \rightarrow \mathsf{E}$ is continuous (see VII.7 of \cite{MM2012}) with respect to the Gaeta topology over $\mathsf{C}$ if and only $G(0)\cong 0$ and $G(1+1)\cong 1+1$; i.e. it preserves binary coproducts.  
\end{remark}

If $\V$ is a coextensive variety, we write $\fpMod(\mathcal{V})$ for the full subcategory of finitely presented algebras of $\V$. As result of Proposition 3.6 \cite{MR2018} and Theorem \ref{fine theorem} we obtain the following:

\begin{proposition}\label{modfp coextensive}
    $\fpMod(\mathcal{V})$ is coextensive.
\end{proposition}

Let $\mathsf{E}$ be a topos. Due to Lawvere's duality \cite{L1963}, it is known that the category of $\V$-models in $\mathsf{E}$ is equivalent to the category of limit preserving functors $\mathsf{Lex}(\fpMod(\mathcal{V})^{\mathsf{op}}, \mathsf{E})$. So, for every $\V$-model $M$ in $\mathsf{E}$, there exists an essentially unique limit preserving functor $\phi_{M}:\fpMod(\mathcal{V})^{\mathsf{op}} \rightarrow \mathsf{E}$, such that $\phi_{M}(\textbf{F}_{\V}(x))\cong M$. 
\\

Now we focus our discussion on $(\vec{0},\vec{1})$-dense coextensive varieties. Let $\V$ be one such variety. Let $x_{1},\ldots,x_{k}$ be a finite set of variables and let $p_{1},\ldots,p_{k},$ $q_{1},\ldots, q_{k}$ be terms in the language of $\V$ with variables $y_{1},\ldots,y_{l}$. If $\delta$ denotes the congruence $\bigvee_{i=1}^{k}\mathsf{Cg}^{\mathbf{F}_{\V}(\vec{y})}(p_{i}(\vec{y}),q_{i}(\vec{y}))$ and $\mathbf{A}$ denotes the algebra $\mathbf{F}_{\V}(\vec{y})/\delta$, observe that for a $\V$-model $M$ in $\mathsf{E}$, the functor $\phi_{M}$ sends the finitely presentable algebra $\mathbf{A}$ to the following equalizer in $\mathsf{E}$:
\[\xymatrix{
\phi_{M}(\mathbf{A}) \ar@{>->}[r] & M^{l} \ar@< 2pt>[rr]^-{\langle p_{M_1},\ldots,p_{M_k}\rangle}
\ar@<-2pt>[rr]_-{\langle q_{M_1},\ldots,q_{M_k}\rangle} & & M^{k}
} \]
where $p_{M_i}$ and $q_{M_i}$ denote the interpretation in $\mathsf{E}$ of the terms $p_{i}$ and $q_{i}$ in $M$, respectively, with $1\leq i\leq k$. I.e. the image of $\mathbf{A}$ by $\phi_{M}$ essentially coincides in $\mathsf{E}$ with the interpretation in $M$ of the formula 
\[\varepsilon(\vec{y}):=\bigwedge_{i=1}^{k} p_{i}(\vec{y})=q_{i}(\vec{y}).\] 
In what follows, we write $\mathcal{G}(\mathcal{V})$ for the Gaeta topos determined by the extensive category $\fpMod(\mathcal{V})^{\mathsf{op}}$. 
\\

Let $M$ be a $\V$-model in a topos $\mathsf{E}$. As result of the above discussion, now we can restate Lemma \ref{Connected VModels in a topos} by means of the functor $\phi_{M}$. 

\begin{lemma}\label{V-indecomposable in a topos}
Let $\V$ be a $(\vec{0},\vec{1})$-dense coextensive variety, $\mathsf{E}$ be a topos and let $M$ be a $\V$-model in $\mathsf{E}$.  Then, the following are equivalent:
\begin{itemize}
\item[(1)] $M$ is central-free in $\mathsf{E}$.
\item[(2)] $\phi_{M}(\mathbf{F}_{\V}(\vec{x},\vec{y})/\theta)\cong 1+1$ and $\phi_{M}(\mathbf{1})\cong 0$.
\end{itemize} 
\end{lemma}
\begin{proof}
Let $M$ be a $\V$-model in a topos $\mathsf{E}$. Notice that by Corollary \ref{decomposition of 0x0}, it is the case that \[[\sigma(\vec{x},\vec{y})]_{M}\cong \phi_{M}(\mathbf{F}_{\V}(\vec{x},\vec{y})/\theta)\cong \phi_{M}(\mathbf{0}\times \mathbf{0}),\] where $\theta$ is defined as in (\ref{definincion teta}), and \[[\vec{0}=\vec{1}]_{M} \cong \phi_{M}(\mathbf{0}/\mathsf{Cg}^{\mathbf{0}}(\vec{0},\vec{1})).\]
Since $\V$ is a variety with $\vec{0}$ and $\vec{1}$, then \[\mathbf{0}/\mathsf{Cg}^{\mathbf{0}}(\vec{0},\vec{1})\cong \mathbf{1}.\]
Hence from Lemma \ref{Connected VModels in a topos} and Remark \ref{Gaeta continous} it is immediate that a $\V$-model $M$ in $\mathsf{E}$ is central-free in such a topos if and only if $(2)$ holds.

\end{proof}

\begin{corollary}
    Let $\mathsf{E}$ be a cocomplete topos and let $\mathsf{D}$ be the full subcategory of  $\mathsf{Lex}(\fpMod(\mathcal{V})^{\mathsf{op}},\mathsf{E})$ whose objects are functors that send $\mathbf{F}_{\V}(\vec{x},\vec{y})/\theta$ in $1+1$ and $\mathbf{1}$ in $0$. Then, there is an equivalence between $\mathsf{Mod}(\mathcal{V}_{CF},\mathsf{E})$ and $\mathsf{D}$.
\end{corollary}

The following Lemma establishes that in a coextensive variety $\mathcal{V}$, finite product decompositions of algebras in $\mathcal{V}$ are determined by ``copartitions" of $\vec{0}$ by central elements. Moreover, as will be seen later, this result, leads to a dual description of the basis for the Gaeta topology in $\mathsf{Mod_{fp}}(\mathcal{V})$, when $\mathcal{V}$ is $(\vec{0},\vec{1})$-dense. It is important to note that the formulation of this result involves the Boolean algebra of central elements, as presented in subsection \ref{Central Elements}. 

\begin{lemma}\label{Description Gaeta basis}
  Let $\V$ be a coextensive variety and let $\mathbf{A}\in \V$. If $\vec{e}_1,...,\vec{e}_n \in Z(\mathbf{A})$ satisfy the following conditions:
  \begin{enumerate}
      \item $\vec{e}_1\wedge_{\mathbf{A}}\ldots \wedge_{\mathbf{A}}\vec{e}_n=\vec{0}$, 
      \item $\vec{e}_i\vee_{\mathbf{A}}\vec{e}_j=\vec{1}$, for every $i\neq j$,
  \end{enumerate}
  Then $\mathbf{A}$ is isomorphic to $\underset{1 \leq i \leq n}{\prod}\mathbf{A}/\mathsf{Cg}^{\mathbf{A}}(\vec{e}_i,\vec{0})$.
\end{lemma}
\begin{proof}
    We proceed by induction on $n$. If $n=2$, the result follows due to well known facts on central elements. If $n=3$, let us assume $\vec{e}_1\wedge_{\mathbf{A}}\vec{e}_2 \wedge_{\mathbf{A}}\vec{e}_3=\vec{0}$ and $\vec{e}_i\vee_{\mathbf{A}}\vec{e}_j=\vec{1}$, for every $i\neq j$. Let $\vec{e}=\vec{e}_1\wedge_{\mathbf{A}}\vec{e}_2$. Straightforward calculations show that $\vec{e}\diamond_{\mathbf{A}} \vec{e}_3$ and therefore, we get that $\mathbf{A}\cong \mathbf{A}/\mathsf{Cg}^{\mathbf{A}}(\vec{e},\vec{0})\times \mathbf{A}/\mathsf{Cg}^{\mathbf{A}}(\vec{e}_3,\vec{0})$. Now observe that Corollary 4 of \cite{V1999} allows to show that $\mathsf{Cg}^{\mathbf{A}}(\vec{e},\vec{0})\subseteq \mathsf{Cg}^{\mathbf{A}}(\vec{e}_1,\vec{0}), \mathsf{Cg}^{\mathbf{A}}(\vec{e}_2,\vec{0})$, $\mathsf{Cg}^{\mathbf{A}}(\vec{e},\vec{0})=\mathsf{Cg}^{\mathbf{A}}(\vec{e}_1,\vec{0})\cap \mathsf{Cg}^{\mathbf{A}}(\vec{e}_2,\vec{0})$ and $\mathsf{Cg}^{\mathbf{A}}(\vec{e}_1,\vec{0})\vee \mathsf{Cg}^{\mathbf{A}}(\vec{e}_1,\vec{0})= \nabla^{\mathbf{A}}$. Hence, by Lemma \ref{factor congruence quotinents}, \[\mathbf{A}/\mathsf{Cg}^{\mathbf{A}}(\vec{e},\vec{0}) \cong \mathbf{A}/\mathsf{Cg}^{\mathbf{A}}(\vec{e}_1,\vec{0}) \times \mathbf{A}/\mathsf{Cg}^{\mathbf{A}}(\vec{e}_2,\vec{0}),\]
    so the statement for $n=3$ holds. The details of the inductive step are left to the reader.
\end{proof}

Now, we are ready to show the main result of this section.

\begin{theorem}
Let $\mathcal{V}$ be a $(\vec{0},\vec{1})$-dense coextensive variety. Then $\mathcal{G}(\mathcal{V})$ is a classifying topos for central-free $\mathcal{V}$-models.   \end{theorem}
\begin{proof}
    We start by stressing that thanks to Lemma \ref{Description Gaeta basis} we are able to provide a dual description of the basis for the Gaeta topology in $\mathsf{Mod}_{\mathsf{fp}}(\V)$. I.e. a Gaeta cocover on a finitely presented algebra of $\V$ is a finite family $\{\mathbf{A}\rightarrow \mathbf{A}/\mathsf{Cg}^{\mathbf{A}}(\vec{e}_i,\vec{0})\mid i\in I\}$ of maps in $\mathsf{Mod}_{\mathsf{fp}}(\V)$ such that (1) and (2) of Lemma \ref{Description Gaeta basis} hold. We recall that, there exists a somewhat standard method for presenting
a specific site for the classifier of central-free $\V$-models. For instance,
Proposition D3.1.10 of \cite{J2002} provides an illustration of this approach. The method consists of two steps: (1) Construct the classifier for the theory presented by the equations of $\mathcal{V}$. In this case, this is the topos of functors $\mathsf{Mod}_{\mathsf{fp}}(\V) \to \Set$. (2) The remaining axioms are 'forced' by introducing an appropriate Grothendieck topology.
 Observe that from Corollary \ref{decomposition of 0x0}, we are able to consider the least Grothendieck topology ``containing" the empty cocover on the terminal object and the cocover 
\begin{displaymath}
    \xymatrix{
   \mathbf{F}_{\V}(\vec{x},\vec{y})/\mu & \mathbf{F}_{\V}(\vec{x},\vec{y})/\theta \ar[l] \ar[r] & \mathbf{F}_{\V}(\vec{x},\vec{y})/\lambda 
    }
\end{displaymath}
with $\theta$, $\mu$ and $\lambda$ as defined in (\ref{definincion teta}). Let us denote by $G$ such a Grothendieck topology. The concrete dual description of the Gaeta topology that we provided implies that the two cocovers generating 
$G$ belong to the basis of the Gaeta topology. Conversely, any binary cocover 
\begin{displaymath}
\xymatrix{
\mathbf{A}/\mathsf{Cg}^{\mathbf{A}}(\vec{e}_1,\vec{0}) & \mathbf{A} \ar[l] \ar[r] & \mathbf{A}/\mathsf{Cg}^{\mathbf{A}}(\vec{e}_2,\vec{0})
}
\end{displaymath}
leads to $\vec{e}_1 \diamond_{\mathbf{A}} \vec{e}_2$, so from the coextensivity of $\mathsf{Mod}_{\mathsf{fp}}(\V)$, it follows that such a cocover arises as a pushout of the cocovers that generate $G$.  The assertion that all non-empty Gaeta cocovers are in $G$ follows from a straightforward inductive argument as outlined in Lemma VIII.6.2 of \cite{MM2012}. Thus, the Gaeta
topology is included in $G$. Therefore, the two topologies coincide. This concludes the proof.
\end{proof}

\section*{Acknowledgements}
The author wishes to convey his appreciation for the institutional support received from the National Scientific and Technical Research Council (CONICET). Furthermore, this project has been funded by the MOSAIC Project 101007627, within the European Union's Horizon 2020 research and innovation programme under the Marie Skłodowska-Curie Actions.

\vline
\\
William Zuluaga,\\
Facultad de Ciencias Exactas (UNCPBA),\\
Pinto 399, Tandil (7000),\\
and CONICET, Argentina,\\
wizubo@gmail.com\\
ORCID: \url{https://orcid.org/0000-0002-8798-9493}

\end{document}